\documentclass[10pt]{article}
\usepackage{cite}
\usepackage{mathrsfs}
\usepackage{amsfonts}
\usepackage{amsmath}
\usepackage{amsfonts,amssymb}
\usepackage{dsfont}
\usepackage{curves}
\usepackage{mathrsfs}
\usepackage{pifont}
\usepackage{amssymb}
\allowdisplaybreaks

\numberwithin{equation}{section}

\date{}

\def\BigRoman{\uppercase\expandafter{\romannumeral\number\count 255 }}
\def\Romannumeral{\afterassignment\BigRoman\count255=}

\begin{document}
\title{Sufficient conditions for a special factor in a graph with minimum degree
}
\author{\small  Sizhong Zhou\footnote{Corresponding author. E-mail address: zsz\_cumt@163.com (S. Zhou)}\\
\small School of Science, Jiangsu University of Science and Technology,\\
\small Zhenjiang, Jiangsu 212100, China\\
}

\maketitle
\begin{abstract}
\noindent Let $G$ be a graph. The size and the signless Laplacian spectral radius of $G$ are denoted by $e(G)$ and $q(G)$, respectively. A spanning subgraph $F$ of $G$ is called an $H_b$-factor of $G$
if $d_F(v)\in\{1,3,5,\ldots,b-1,b\}$ for every $v\in V(G)$, where $b\geq2$ is an even integer. Lu and Wang obtained a sufficient condition according to the number of odd components in $G-S$ for a connected
graph $G$ of even order to have an $H_b$-factor, where $S$ is a subset of $V(G)$ [H. Lu, D. Wang, On Cui-Kano's characterization problem on graph factors, J. Graph Theory 74 (2013) 335--343]. In this paper,
motivated by Lu and Wang's above result, we establish a lower bound for the size in an $n$-vertex connected graph $G$ with given minimum degree to guarantee that $G$ has an $H_b$-factor. Further, we show
a lower bound for the signless Laplacian spectral radius in an $n$-vertex 2-connected graph $G$ with given minimum degree to ensure that $G$ has an $H_b$-factor.
\\
\begin{flushleft}
{\em Keywords:} graph; minimum degree; size; signless Laplacian spectral radius; $H_b$-factor.

(2020) Mathematics Subject Classification: 05C70, 05C50
\end{flushleft}
\end{abstract}

\section{Introduction}

Let $G=(V(G),E(G))$ denote a finite, undirected and simple graph, where $V(G)$ denotes the vertex set of $G$ and $E(G)$ denotes the edge set of $G$. Let $|E(G)|=e(G)$ and $|V(G)|=n$ denote the size and the order
of $G$, respectively. For any $v\in V(G)$, let $d_G(v)$ denote the degree of $v$ in $G$. We denote by $\delta(G)$ and $o(G)$ the minimum degree and the number of odd components in $G$, respectively. For
$S\subseteq V(G)$, the subgraph of $G$ induced by $S$ is denoted by $G[S]$, while the subgraph of $G$ induced by $V(G)-S$ is denoted by $G-S$. We denote by $K_n$ the complete graph of order $n$. Let $G_1$ and
$G_2$ be two vertex-disjoint graphs. Then the union of $G_1$ and $G_2$ is denoted by $G_1\cup G_2$. The join $G_1\vee G_2$ is the graph obtained from $G_1\cup G_2$ by adding every possible edge between $V(G_1)$
and $V(G_2)$.

Let $A(G)$ and $D(G)$ denote the adjacency matrix and the diagonal degree matrix of $G$, respectively. The signless Laplacian matrix of $G$, denoted by $Q(G)$, is defined by $Q(G)=D(G)+A(G)$. The largest eigenvalue
of $Q(G)$ is called the signless Laplacian spectral radius of $G$, which is denoted by $q(G)$. We refer reader to \cite{BOT,Oa,Ws,ZLF,ZBS,ZZZL,WZL} for some properties on the signless Laplacian spectral radius
of $G$.

Let $b\geq2$ be an even integer. A spanning subgraph $F$ of $G$ is called an odd $[1,b-1]$-factor of $G$ if $d_F(v)\in\{1,3,5,\ldots,b-1\}$ for each $v\in V(G)$. An odd $[1,b-1]$-factor of $G$ is called a perfect
matching of $G$ if $b=2$. A spanning subgraph $F$ of $G$ is called an $H_b$-factor of $G$ if $d_F(v)\in\{1,3,5,\ldots,b-1,b\}$ for every $v\in V(G)$.

O \cite{Os} proposed two sufficient conditions for a graph with a perfect matching based on the size and the adjacency spectral radius. Fan, Lin and Lu \cite{FLL}, Zhou and Liu \cite{ZL} obtained some results on
the existence of odd $[1,b-1]$-factors in graphs. Zhou \cite{Zs3} established sharp lower bounds for both the size and the adjacency spectral radius in a graph $G$ with given minimum degree to ensure that $G-W$
contains an odd $[1,b-1]$-factor for any $W\subseteq V(G)$ with $|W|=k$. Wang, Yang and Yang \cite{WYY} presented a sufficient condition involving the signless Laplacian spectral radius to ensure that $G-W$ has
an odd $[1,b-1]$-factor for any $W\subseteq V(G)$ with $|W|=k$. Lu and Wang \cite{LW} provided some sufficient conditions for a graph to possess an $H_b$-factor. Fan and Liu \cite{FL} claimed an adjacency spectral
radius condition for a graph with an $H_b$-factor. We refer reader to \cite{KT,Wc,Wu,PZ,Zs1,Zs2,ZBW,Zr,Zhou,ZZL,Oe,KO} for some other results on graph factors.

Motivated by \cite{LW} directly, we put forward two sufficient conditions by virtue of the size and the signless Laplacian spectral radius to guarantee that a graph $G$ of even order $n$ with minimum degree $\delta$
has an $H_b$-factor.

\medskip

\noindent{\textbf{Theorem 1.1.}} Let $G$ be a connected graph of even order $n\geq\max\{\frac{b^{2}\delta+2b\delta-\delta-1}{b-1},\frac{b^{2}\delta^{2}+(4b^{2}+5b)\delta+b^{2}+4b+1}{6b}\}$ with minimum degree
$\delta\geq1$, where $b\geq2$ is an even integer. If
$$
e(G)\geq e(K_{\delta}\vee(K_{n-(b+1)\delta}\cup b\delta K_1)),
$$
then $G$ has an $H_b$-factor, unless $G=K_{\delta}\vee(K_{n-(b+1)\delta}\cup b\delta K_1)$.

\medskip

\noindent{\textbf{Theorem 1.2.}} Let $G$ be a 2-connected graph of even order $n\geq\max\{4b\delta,b\delta^{3}-b\delta-2\}$ with minimum degree $\delta\geq2$, where $b\geq2$ is an even integer. If
$$
q(G)\geq q(K_{\delta}\vee(K_{n-(b+1)\delta}\cup b\delta K_1)),
$$
then $G$ has an $H_b$-factor, unless $G=K_{\delta}\vee(K_{n-(b+1)\delta}\cup b\delta K_1)$.

\medskip

\section{Some preliminaries}

In 2013, Lu and Wang \cite{LW} posed a sufficient condition for the existence of an $H_b$-factor in a graph with even order.

\medskip

\noindent{\textbf{Lemma 2.1}} (Lu and Wang \cite{LW}). Let $b\geq2$ be an even integer, and let $G$ be a connected graph of even order. If
$$
o(G-S)\leq b|S|
$$
for every nonempty subset $S\subseteq V(G)$, then $G$ has an $H_b$-factor.

\medskip

\noindent{\textbf{Lemma 2.2}} (Zheng, Li, Luo and Wang \cite{ZLLW}). Let $\sum\limits_{i=1}^{t}n_i=n-s$. If $n_1\geq n_2\geq\cdots\geq n_t\geq p\geq1$, then
$$
e(K_s\vee(K_{n_1}\cup K_{n_2}\cup\cdots\cup K_{n_t}))\leq e(K_s\vee(K_{n-s-p(t-1)}\cup(t-1)K_p)),
$$
with equality following if and only if $(n_1,n_2,\ldots,n_t)=(n-s-p(t-1),p,\ldots,p)$.

\medskip

\noindent{\textbf{Lemma 2.3}} (Shen, You, Zhang and Li \cite{SYZL}). Let $H$ be a subgraph of a connected graph $G$. Then
$$
q(G)\geq q(H),
$$
with equality occurring if and only if $G=H$.

\medskip

\noindent{\textbf{Lemma 2.4}} (Zheng, Li, Luo and Wang \cite{ZLLW}). Let $\sum\limits_{i=1}^{t}n_i=n-s$. If $n_1\geq n_2\geq\cdots\geq n_t\geq p\geq1$, then
$$
q(K_s\vee(K_{n_1}\cup K_{n_2}\cup\cdots\cup K_{n_t}))\leq q(K_s\vee(K_{n-s-p(t-1)}\cup(t-1)K_p)),
$$
with equality holding if and only if $(n_1,n_2,\ldots,n_t)=(n-s-p(t-1),p,\ldots,p)$.

\medskip

Let $\mathcal{N}=\{1,2,\ldots,n\}$. For a real $n\times n$ matrix $M$ and a partition $\pi:\mathcal{N}=\mathcal{N}_1\cup\mathcal{N}_2\cup\cdots\cup\mathcal{N}_r$, the matrix $M$ can be written as
\begin{align*}
M=\left(
  \begin{array}{cccc}
    M_{11} & M_{12} & \cdots & M_{1r}\\
    M_{21} & M_{22} & \cdots & M_{2r}\\
    \vdots & \vdots & \ddots & \vdots\\
    M_{r1} & M_{r2} & \cdots & M_{rr}\\
  \end{array}
\right).
\end{align*}
The quotient matrix, denoted by $M_{\pi}$, corresponding to the partition $\pi$ is defined as $M_{\pi}=(m_{ij})_{r\times r}$, where $m_{ij}$ is the average row sum of $M_{ij}$. The partition $\pi$ is equitable if
the row sum of every $M_{ij}$ is a constant for $1\leq i,j\leq r$.

\medskip

\noindent{\textbf{Lemma 2.5}} (You, Yang, So and Xi\cite{YYSX}). Let $M$ be a real $n\times n$ matrix with an equitable partition $\pi$, and let $M_{\pi}$ be the corresponding quotient matrix. Then every eigenvalue
of $M_{\pi}$ is an eigenvalue of $M$. Furthermore, if $M$ is nonnegative and irreducible, then the largest eigenvalues of $M$ and $M_{\pi}$ are equal.

\section{The proof of Theorem 1.1}

\noindent{\it Proof of Theorem 1.1.} Suppose that $G$ has no $H_b$-factor. By Lemma 2.1, there exists a nonempty subset $S\subseteq V(G)$ satisfying $o(G-S)\geq b|S|+1$. Let $|S|=s$. It is obvious that $G$ is a
spanning subgraph of $G_1=K_s\vee(K_{n_1}\cup K_{n_2}\cup\cdots\cup K_{n_{bs+1}})$, where $n_1\geq n_2\geq\cdots\geq n_{bs+1}$ are positive odd integers with $\sum\limits_{i=1}^{bs+1}n_i=n-s$. Then the following
inequality holds:
\begin{align}\label{eq:3.1}
e(G)\leq e(G_1),
\end{align}
with equality occurring if and only if $G=G_1$. The rest of the proof is carried out by analyzing three different cases based on the value of $s$.

\noindent{\bf Case 1.} $s\geq\delta+1$.

Let $G_2=K_s\vee(K_{n-(b+1)s}\cup bsK_1)$, where $n\geq(b+1)s+1$. According to Lemma 2.2, we infer
\begin{align}\label{eq:3.2}
e(G_1)\leq e(G_2),
\end{align}
with equality holding if and only if $(n_1,n_2,\ldots,n_{bs+1})=(n-(b+1)s,1,\ldots,1)$.

Define $G_*=K_{\delta}\vee(K_{n-(b+1)\delta}\cup b\delta K_1)$. It follows from $b\geq2$, $s\geq\delta+1$, $n\geq(b+1)s+1$ and $n\geq\frac{b^{2}\delta+2b\delta-\delta-1}{b-1}$ that
\begin{align*}
e(G_*)-e(G_2)=&\binom{n-b\delta}{2}+b\delta^{2}-\binom{n-bs}{2}-bs^{2}\\
=&\frac{1}{2}(s-\delta)(2bn-b^{2}s-2bs-b^{2}\delta-2b\delta-b)\\
=&\frac{1}{2}(s-\delta)((b-1)n+(b+1)n-b^{2}s-2bs-b^{2}\delta-2b\delta-b)\\
\geq&\frac{1}{2}(s-\delta)((b-1)n+(b+1)((b+1)s+1)-b^{2}s-2bs-b^{2}\delta-2b\delta-b)\\
=&\frac{1}{2}(s-\delta)((b-1)n+s-b^{2}\delta-2b\delta+1)\\
\geq&\frac{1}{2}(s-\delta)(b^{2}\delta+2b\delta-\delta-1+\delta+1-b^{2}\delta-2b\delta+1)\\
=&\frac{1}{2}(s-\delta)\\
>&0,
\end{align*}
which leads to $e(G_2)<e(G_*)$. Combining this with \eqref{eq:3.1} and \eqref{eq:3.2}, we have
$$
e(G)\leq e(G_1)\leq e(G_2)<e(G_*)=e(K_{\delta}\vee(K_{n-(b+1)\delta}\cup b\delta K_1)),
$$
which contradicts $e(G)\geq e(K_{\delta}\vee(K_{n-(b+1)\delta}\cup b\delta K_1))$.

\noindent{\bf Case 2.} $s=\delta$.

In this case, $G_1=K_{\delta}\vee(K_{n_1}\cup K_{n_2}\cup\cdots\cup K_{n_{b\delta+1}})$ for some positive odd integers $n_1\geq n_2\geq\cdots\geq n_{b\delta+1}$ with $\sum\limits_{i=1}^{b\delta+1}n_i=n-\delta$.
By means of Lemma 2.2, we conclude
$$
e(G_1)\leq e(K_{\delta}\vee(K_{n-(b+1)\delta}\cup b\delta K_1)),
$$
with equality if and only if $(n_1,n_2,\ldots,n_{b\delta+1})=(n-(b+1)\delta,1,\ldots,1)$. Combining this with \eqref{eq:3.1}, we obtain
$$
e(G)\leq e(K_{\delta}\vee(K_{n-(b+1)\delta}\cup b\delta K_1)),
$$
with equality if and only if $G=K_{\delta}\vee(K_{n-(b+1)\delta}\cup b\delta K_1)$, a contradiction.

\noindent{\bf Case 3.} $s\leq\delta-1$.

Define $G_3=K_s\vee(K_{n-s-bs(\delta+1-s)}\cup bsK_{\delta+1-s})$. Recall that $G$ is a spanning subgraph of $G_1=K_s\vee(K_{n_1}\cup K_{n_2}\cup\cdots\cup K_{n_{bs+1}})$ for positive odd integers
$n_1\geq n_2\geq\cdots\geq n_{bs+1}$ with $\sum\limits_{i=1}^{bs+1}n_i=n-s$. Since $\delta(G_1)\geq\delta(G)=\delta$, we deduce $n_{bs+1}\geq\delta+1-s$. In terms of Lemma 2.2, we get
\begin{align}\label{eq:3.3}
e(G_1)\leq e(G_3),
\end{align}
where the equality occurs if and only if $(n_1,n_2,\ldots,n_{bs+1})=(n-s-bs(\delta+1-s),\delta+1-s,\ldots,\delta+1-s)$.

Recall that $G_*=K_{\delta}\vee(K_{n-(b+1)\delta}\cup b\delta K_1)$. Then we possess
\begin{align}\label{eq:3.4}
e(G_*)-e(G_3)=&\binom{n-b\delta}{2}+b\delta^{2}-\binom{n-bs(\delta+1-s)}{2}-bs\binom{\delta+1-s}{2}-bs^{2}(\delta+1-s)\nonumber\\
=&\frac{1}{2}(\delta-s)(b^{2}s^{3}-(b^{2}\delta+2b^{2}+b)s^{2}+(2bn-b\delta+b^{2})s-2bn+b^{2}\delta+2b\delta+b).
\end{align}
Let $\varphi(x)=b^{2}x^{3}-(b^{2}\delta+2b^{2}+b)x^{2}+(2bn-b\delta+b^{2})x-2bn+b^{2}\delta+2b\delta+b$. Then the derivative function of $\varphi(x)$ is given by
$$
\varphi'(x)=3b^{2}x^{2}-2(b^{2}\delta+2b^{2}+b)x+2bn-b\delta+b^{2}.
$$
Obviously, the symmetry axis of $\varphi'(x)$ is $x=\frac{b\delta+2b+1}{3b}$. Thus, we deduce $\varphi'(x)\geq\varphi'(\frac{b\delta+2b+1}{3b})=\frac{1}{3}(6bn-b^{2}\delta^{2}-(4b^{2}+5b)\delta-b^{2}-4b-1)\geq0$
due to $n\geq\frac{b^{2}\delta^{2}+(4b^{2}+5b)\delta+b^{2}+4b+1}{6b}$. This implies that $\varphi(x)$ is increasing in the interval $[1,\delta-1]$. Together with $1\leq s\leq\delta-1$, we deduce
\begin{align}\label{eq:3.5}
\varphi(s)\geq\varphi(1)=b\delta>0.
\end{align}
By virtue of \eqref{eq:3.4}, \eqref{eq:3.5} and $s\leq\delta-1$, we infer
\begin{align}\label{eq:3.6}
e(G_3)<e(G_*).
\end{align}
It follows from \eqref{eq:3.1}, \eqref{eq:3.3} and \eqref{eq:3.6} that
$$
e(G)\leq e(G_1)\leq e(G_3)<e(G_*)=e(K_{\delta}\vee(K_{n-(b+1)\delta}\cup b\delta K_1)),
$$
which is a contradiction to $e(G)\geq e(K_{\delta}\vee(K_{n-(b+1)\delta}\cup b\delta K_1))$. This completes the proof of Theorem 1.1. \hfill $\Box$

\medskip

\section{The proof of Theorem 1.2}

\noindent{\it Proof of Theorem 1.2.} Suppose that $G$ has no $H_b$-factor. Based on Lemma 2.1, there exists a nonempty subset $S\subseteq V(G)$ satisfying $o(G-S)\geq b|S|+1$. Let $|S|=s$. Then $G$ is a spanning
subgraph of $G_1=K_s\vee(K_{n_1}\cup K_{n_2}\cup\cdots\cup K_{n_{bs+1}})$ for some positive odd integers $n_1\geq n_2\geq\cdots\geq n_{bs+1}$ with $\sum\limits_{i=1}^{bs+1}n_i=n-s$. In view of Lemma 2.3, we possess
\begin{align}\label{eq:4.1}
q(G)\leq q(G_1),
\end{align}
where the equality follows if and only if $G=G_1$. Next, we shall consider three cases according to the value of $s$.

\noindent{\bf Case 1.} $s\geq\delta+1$.

Define $G_2=K_s\vee(K_{n-(b+1)s}\cup bsK_1)$, where $n\geq(b+1)s+1$. Using Lemma 2.4, we conclude
\begin{align}\label{eq:4.2}
q(G_1)\leq q(G_2),
\end{align}
with equality if and only if $(n_1,n_2,\ldots,n_{bs+1})=(n-(b+1)s,1,\ldots,1)$. The quotient matrix of $Q(G_2)$ in terms of the partition $V(G_2)=V(K_s)\cup V(K_{n-(b+1)s})\cup V(bsK_1)$ is
\begin{align*}
B_2=\left(
  \begin{array}{ccc}
    n+s-2 & n-(b+1)s & bs\\
    s & 2n-(2b+1)s-2 & 0\\
    s & 0 & s\\
  \end{array}
\right),
\end{align*}
and the characteristic polynomial of $B_2$ is
\begin{align*}
f_{B_2}(x)=&x^{3}+(-3n+(2b-1)s+4)x^{2}+(2n^{2}-((2b-3)s+6)n-4bs^{2}+4bs-4s+4)x\\
&-2sn^{2}+(4bs^{2}+6s)n-2b^{2}s^{3}-6bs^{2}-4s.
\end{align*}

Define $G_*=K_{\delta}\vee(K_{n-(b+1)\delta}\cup b\delta K_1)$. Then we denote the quotient matrix of $Q(G_*)$ based on the partition $V(G_*)=V(K_{\delta})\cup V(K_{n-(b+1)\delta})\cup V(b\delta K_1)$ by $B_*$. The
characteristic polynomial of $B_*$ is
\begin{align*}
f_{B_*}(x)=&x^{3}+(-3n+(2b-1)\delta+4)x^{2}+(2n^{2}-((2b-3)\delta+6)n-4b\delta^{2}+4b\delta-4\delta+4)x\\
&-2\delta n^{2}+(4b\delta^{2}+6\delta)n-2b^{2}\delta^{3}-6b\delta^{2}-4\delta.
\end{align*}
Hence, we obtain
\begin{align}\label{eq:4.3}
f_{B_2}(x)-f_{B_*}(x)=(s-\delta)g(x),
\end{align}
where $g(x)=(2b-1)x^{2}-((2b-3)n+4bs+4b\delta-4b+4)x-2n^{2}+(4bs+4b\delta+6)n-2b^{2}s^{2}-(2b^{2}\delta+6b)s-2b^{2}\delta^{2}-6b\delta-4$. The symmetry axis of $g(x)$ is $x=\frac{(2b-3)n+4bs+4b\delta-4b+4}{2(2b-1)}$.
Notice that
$$
\frac{(2b-3)n+4bs+4b\delta-4b+4}{2(2b-1)}<2(n-b\delta-1)
$$
due to $s\geq\delta+1$, $n\geq(b+1)s+1$ and $n\geq4b\delta$. This implies that $g(x)$ is increasing for $x\geq2(n-b\delta-1)$. When $x\geq2(n-b\delta-1)$, it follows from $\delta>1$, $b\geq2$, $s\leq\frac{n-1}{b+1}$
and $n\geq4b\delta$ that
\begin{align}\label{eq:4.4}
g(x)\geq&g(2(n-b\delta-1))\nonumber\\
=&(2b-1)(2(n-b\delta-1))^{2}-((2b-3)n+4bs+4b\delta-4b+4)(2(n-b\delta-1))\nonumber\\
&-2n^{2}+(4bs+4b\delta+6)n-2b^{2}s^{2}-(2b^{2}\delta+6b)s-2b^{2}\delta^{2}-6b\delta-4\nonumber\\
=&4bn^{2}-(4bs+12b^{2}\delta+2b\delta+4b)n-2b^{2}s^{2}\nonumber\\
&+(6b^{2}\delta+2b)s+8b^{3}\delta^{2}+2b^{2}\delta^{2}+8b^{2}\delta+2b\delta\nonumber\\
\geq&4bn^{2}-\Big(4b\Big(\frac{n-1}{b+1}\Big)+12b^{2}\delta+2b\delta+4b\Big)n-2b^{2}\Big(\frac{n-1}{b+1}\Big)^{2}\nonumber\\
&+(6b^{2}\delta+2b)\Big(\frac{n-1}{b+1}\Big)+8b^{3}\delta^{2}+2b^{2}\delta^{2}+8b^{2}\delta+2b\delta\nonumber\\
=&\frac{1}{(b+1)^{2}}\Big((4b^{3}+2b^{2})n^{2}-(12b^{4}\delta+20b^{3}\delta+10b^{2}\delta+2b\delta+4b^{3}-2b^{2}-2b)n\nonumber\\
&+8b^{5}\delta^{2}+18b^{4}\delta^{2}+12b^{3}\delta^{2}+2b^{2}\delta^{2}+8b^{4}\delta+12b^{3}\delta+6b^{2}\delta+2b\delta-4b^{2}-2b\Big)\nonumber\\
\geq&\frac{1}{(b+1)^{2}}\Big((4b^{3}+2b^{2})(4b\delta)^{2}-(12b^{4}\delta+20b^{3}\delta+10b^{2}\delta+2b\delta+4b^{3}-2b^{2}-2b)(4b\delta)\nonumber\\
&+8b^{5}\delta^{2}+18b^{4}\delta^{2}+12b^{3}\delta^{2}+2b^{2}\delta^{2}+8b^{4}\delta+12b^{3}\delta+6b^{2}\delta+2b\delta-4b^{2}-2b\Big)\nonumber\\
=&\frac{1}{(b+1)^{2}}\Big((24b^{5}-30b^{4}-28b^{3}-6b^{2})\delta^{2}+(-8b^{4}+20b^{3}+14b^{2}+2b)\delta-4b^{2}-2b\Big)\nonumber\\
>&\frac{1}{(b+1)^{2}}\Big((24b^{5}-30b^{4}-28b^{3}-6b^{2})\delta+(-8b^{4}+20b^{3}+14b^{2}+2b)\delta-4b^{2}-2b\Big)\nonumber\\
=&\frac{1}{(b+1)^{2}}\Big((24b^{5}-38b^{4}-8b^{3}+8b^{2}+2b)\delta-4b^{2}-2b\Big)\nonumber\\
>&\frac{1}{(b+1)^{2}}(24b^{5}-38b^{4}-8b^{3}+4b^{2})\nonumber\\
>&0.
\end{align}
In terms of \eqref{eq:4.3}, \eqref{eq:4.4} and $s\geq\delta+1$, we obtain $f_{B_2}(x)>f_{B_*}(x)$ for $x\geq2(n-b\delta-1)$. Notice that $K_{n-b\delta}$ is a proper subgraph of $G_*=K_{\delta}\vee(K_{n-(b+1)\delta}\cup b\delta K_1)$.
Together with Lemma 2.3, we have $q(G_*)>q(K_{n-b\delta})=2(n-b\delta-1)$. Thus, we deduce
\begin{align}\label{eq:4.5}
q(G_2)<q(G_*),
\end{align}
where $q(G_2)$ and $q(G_*)$ respectively are the largest roots of $f_{B_2}(x)=0$ and $f_{B_*}(x)=0$ due to Lemma 2.5. From \eqref{eq:4.1}, \eqref{eq:4.2} and \eqref{eq:4.5}, we infer
$$
q(G)\leq q(G_1)\leq q(G_2)<q(G_*)=q(K_{\delta}\vee(K_{n-(b+1)\delta}\cup b\delta K_1)),
$$
which is a contradiction to $q(G)\geq q(K_{\delta}\vee(K_{n-(b+1)\delta}\cup b\delta K_1))$.

\noindent{\bf Case 2.} $s=\delta$.

Recall that $G_1=K_s\vee(K_{n_1}\cup K_{n_2}\cup\cdots\cup K_{n_{bs+1}})$. Using $s=\delta$ and Lemma 2.4, we conclude
$$
q(G_1)\leq q(K_{\delta}\vee(K_{n-(b+1)\delta}\cup b\delta K_1)),
$$
with equality if and only if $G_1=K_{\delta}\vee(K_{n-(b+1)\delta}\cup b\delta K_1)$. Together with \eqref{eq:4.1}, we possess
$$
q(G)\leq q(K_{\delta}\vee(K_{n-(b+1)\delta}\cup b\delta K_1)),
$$
where the equality occurs if and only if $G=K_{\delta}\vee(K_{n-(b+1)\delta}\cup b\delta K_1)$, a contradiction.

\noindent{\bf Case 3.} $s\leq\delta-1$.

Define $G_3=K_s\vee(K_{n-s-bs(\delta+1-s)}\cup bsK_{\delta+1-s})$. Recall that $G$ is a spanning subgraph of $G_1=K_s\vee(K_{n_1}\cup K_{n_2}\cup\cdots\cup K_{n_{bs+1}})$ for some positive odd integers
$n_1\geq n_2\geq\cdots\geq n_{bs+1}$ with $\sum\limits_{i=1}^{bs+1}n_i=n-s$. Since $\delta(G_1)\geq\delta(G)=\delta$, we conclude $n_{bs+1}\geq\delta+1-s$. Using Lemma 2.4, we obtain
\begin{align}\label{eq:4.6}
q(G_1)\leq q(G_3),
\end{align}
with equality if and only if $(n_1,n_2,\ldots,n_{bs+1})=(n-s-bs(\delta+1-s),\delta+1-s,\ldots,\delta+1-s)$. Note that $n$ is even and $G$ is 2-connected. If $s=1$, then $1=o(G-S)\geq b|S|+1=b+1>1$, which
is a contradiction. Next, we deal with $2\leq s\leq\delta-1$.

We are to show $q(G_3)=q'<2(n-bs(\delta+1-s))$. Assume that $q'\geq2(n-bs(\delta+1-s))$. The Perron vector of $Q(G_3)$ is written as $x$. Using symmetry, $x$ takes the same value (say $x_1$, $x_2$ and $x_3$)
based on the vertices of $K_s$, $K_{n-s-bs(\delta+1-s)}$ and $bsK_{\delta+1-s}$, respectively. According to $Q(G_3)x=q'x$, we deduce
\begin{align}\label{eq:4.7}
q'x_1=(n+s-2)x_1+(n-s-bs(\delta+1-s))x_2+bs(\delta+1-s)x_3,
\end{align}
\begin{align}\label{eq:4.8}
q'x_2=sx_1+(2n-s-2bs(\delta+1-s)-2)x_2
\end{align}
and
\begin{align}\label{eq:4.9}
q'x_3=sx_1+(2\delta-s)x_3.
\end{align}
In terms of \eqref{eq:4.8} and \eqref{eq:4.9}, we possess
\begin{align}\label{eq:4.10}
x_2=\frac{sx_1}{q'-(2n-s-2bs(\delta+1-s)-2)}
\end{align}
and
\begin{align}\label{eq:4.11}
x_3=\frac{sx_1}{q'-(2\delta-s)}.
\end{align}
Using $\delta\geq s+1\geq3$ and $n\geq b\delta^{3}-b\delta-2$, we possess $q'\geq2(n-bs(\delta+1-s))>2(\delta+1)$. Together with \eqref{eq:4.7}, \eqref{eq:4.10}, \eqref{eq:4.11}, $2\leq s\leq\delta-1$ and
$n\geq b\delta^{3}-b\delta-2$, we infer
\begin{align*}
q'=&n+s-2+\frac{s(n-s-bs(\delta+1-s))}{q'-(2n-s-2bs(\delta+1-s)-2)}+\frac{bs^{2}(\delta+1-s)}{q'-(2\delta-s)}\\
<&n+s-2+\frac{s(n-s-bs(\delta+1-s))}{s+2}+\frac{bs^{2}(\delta+1-s)}{s+2}\\
=&n+s-2+\frac{s(n-s)}{s+2}\\
=&2(n-bs(\delta+1-s))-\frac{2n-2bs(s+2)(\delta+1-s)+4}{s+2}\\
\leq&2(n-bs(\delta+1-s))\\
\leq&q',
\end{align*}
which is a contradiction. Hence, we have
\begin{align}\label{eq:4.12}
q'<2(n-bs(\delta+1-s)).
\end{align}

Recall that $q(G_*)>2(n-b\delta-1)$. In terms of \eqref{eq:4.1}, \eqref{eq:4.6}, \eqref{eq:4.12}, $b\geq2$ and $2\leq s\leq\delta-1$, we possess
\begin{align*}
q(G)\leq&q(G_1)\leq q(G_3)=q'<2(n-bs(\delta+1-s))\\
=&2(n-b\delta-1)-2(b(s-1)(\delta-s)-1)\\
<&2(n-b\delta-1)\\
<&q(K_{\delta}\vee(K_{n-(b+1)\delta}\cup b\delta K_1)),
\end{align*}
which contradicts $q(G)\geq q(K_{\delta}\vee(K_{n-(b+1)\delta}\cup b\delta K_1))$. This completes the proof of Theorem 1.2. \hfill $\Box$

\medskip

\section*{Data availability statement}

My manuscript has no associated data.

\section*{Declaration of competing interest}

The authors declare that they have no conflicts of interest to this work.

\section*{Acknowledgments}

This work was supported by the Natural Science Foundation of Jiangsu Province (Grant No. BK20241949).

\end{document}